\documentclass[a4paper,11pt]{article}
\usepackage{amsmath,amssymb, amsbsy, amsfonts, amsthm}
\usepackage{mathrsfs}
\usepackage{bigints}
\usepackage{hyperref}
\usepackage{color,psfrag}
\usepackage[dvips]{graphicx}
\usepackage{enumerate}
\usepackage{tikz}
\usetikzlibrary{decorations.markings}
\tikzset{
  set arrow inside/.code={\pgfqkeys{/tikz/arrow inside}{#1}},
  set arrow inside={end/.initial=>, opt/.initial=},
  /pgf/decoration/Mark/.style={
    mark/.expanded=at position #1 with
    {
      \noexpand\arrow[\pgfkeysvalueof{/tikz/arrow inside/opt}]{\pgfkeysvalueof{/tikz/arrow inside/end}}
    }
  },
  arrow inside/.style 2 args={
    set arrow inside={#1},
    postaction={
      decorate,decoration={
        markings,Mark/.list={#2}
      }
    }
  },
}
\tikzset{
  set arrow outside/.code={\pgfqkeys{/tikz/arrow outside}{#1}},
  set arrow outside={end/.initial=<, opt/.initial=},
  /pgf/decoration/Mark/.style={
    mark/.expanded=at position #1 with
    {
      \noexpand\arrow[\pgfkeysvalueof{/tikz/arrow outside/opt}]{\pgfkeysvalueof{/tikz/arrow outside/end}}
    }
  },
  arrow outside/.style 2 args={
    set arrow outside={#1},
    postaction={
      decorate,decoration={
        markings,Mark/.list={#2}
      }
    }
  },
}
\usepackage{url}
\usepackage{comment}
\usepackage{colortbl}

\usepackage{graphicx,color}
\graphicspath{{Figures/}}
\usepackage{tikz,pgfplots,pgf}
\usetikzlibrary{patterns}
\usepackage[font=small]{caption}
\usepackage{subcaption}

\newcommand{\eps}{\varepsilon}

\renewcommand{\geq }{\geqslant}

\renewcommand{\leq }{\leqslant}

\def\neweq#1{\begin{equation}\label{#1}}
\def\endeq{\end{equation}}

\newtheorem{theorem}{Theorem}
\newtheorem{proposition}[theorem]{Proposition}
\newtheorem{lemma}[theorem]{Lemma}
\newtheorem{corollary}[theorem]{Corollary}
\newtheorem{definition}[theorem]{Definition}

\titlepage 
\title{The role of convection in the limit shape of the critical front profile \\ for Born-Infeld diffusion models}
\author{Maurizio Garrione\thanks{Dipartimento di Matematica, Politecnico di Milano, Milano, Italy}\,,\  
Mohamed Jleli\footnote{Department of Mathematics, King Saud University, Riyadh, Saudi Arabia}\,, \ Bessem Samet\footnote{Department of Mathematics, King Saud University, Riyadh, Saudi Arabia}}  
\date{}

\begin{document}
\maketitle

\begin{abstract}
In this paper, we deal with models with Born-Infeld type diffusion and monostable reaction, investigating the effect of the introduction of a convection term on the limit shape of the critical front profile for vanishing diffusion. We first provide an estimate of the critical speed and then, through a careful analysis of an equivalent first-order problem, we show that different convection terms may lead either to a complete sharpening of the limit profile or to its complete regularization, presenting some related numerical simulations.
\end{abstract}

\smallbreak

\section{Introduction}

This paper is devoted to the analysis of the 
limit shape, for $\eps \to 0^+$, of the critical front profile for the $1$-dimensional 
reaction-convection-diffusion equation
\begin{equation}\label{MinkIntro}
u_t=\eps \left(\frac{u_x}{\sqrt{1-u_x^2}}\right)_x  - (h(u))_x + f(u), \quad u=u(x, t), \; x \in \mathbb{R}, \, t \in \mathbb{R}.
\end{equation}
In the course, we assume that the reaction term $f$ is monostable, namely it fulfills $f(0)=0=f(1)$ and $f(s) > 0$ for every $s \in (0, 1)$ (see assumption (F) in Section \ref{sez2}). 

The second-order operator appearing in equation \eqref{MinkIntro}, known as the \emph{relativistic} or \emph{Born-Infeld} operator, makes the equation singular and constrains regular solutions to be subject to the a priori bound $\vert u_x \vert < 1$. It models diffusion in the Lorentz-Minkowski space (see, e.g., the references in \cite{Azz06, GarPP}) and may be framed more generally inside the theory of \emph{singular diffusion operators}. At the same time, it appears in Born-Infeld electrodynamics \cite{BorInf}, as mentioned in \cite{GarPP, Kru2}. The common feature of these two viewpoints is the presence of a universal finite bound for the key quantities of the model (the speed in the former case, the electric field in the latter one). 

We are concerned with strictly increasing regular wave fronts $u(x, t)=v(x+ct)$ for \eqref{MinkIntro}, solving the second-order problem 
\begin{equation}\label{IIordMink}
\left\{
\begin{array}{l}
\eps \displaystyle \left(\frac{v'}{\sqrt{1-(v')^2}}\right)' - (c+h'(v)) v' + f(v) = 0 \vspace{0.1cm}\\
v(-\infty)=0, \, v(+\infty)=1, \, v' > 0;
\end{array}
\right.
\end{equation}
we impose since the very beginning the condition $v(0)=1/2$, in order to recover uniqueness (notice that \eqref{IIordMink} is autonomous). 
Writing the differential equation in \eqref{IIordMink} as the equivalent first-order system in the phase plane
$$
\left\{
\begin{array}{l}
v'=\dfrac{w}{\sqrt{1+w^2}} \vspace{0.1cm}\\
\eps w'=(c+h'(v))\dfrac{w}{\sqrt{1+w^2}}  - f(v),
\end{array}
\right.
$$
one has 
$$
\frac{dw}{dv} =\frac{c+h'(v)}{\eps} - f(v) \frac{\sqrt{1+w^2}}{\eps w}, 
$$
so that $y(v)=\eps(\sqrt{1+w(v)^2}-1)$ satisfies the first-order two-point problem 
\begin{equation}\label{notazioney}
\left\{
\begin{array}{l}
\displaystyle y'= (c+h'(v)) \frac{\sqrt{y(2\eps+y)}}{\eps+y} - f(v) \vspace{0.1cm}\\
y(0)=0, \, y(1)= 0, \, y > 0 \text{ on } (0, 1)
\end{array}
\right.
\end{equation}
(see also \cite{FifMcL, GarSan, MalMar03}).
Explicitly,
\begin{equation}\label{costruzionev}
y(v)=\eps \left(\frac{1}{\sqrt{1-v'(z(v))^2}}-1\right) \quad \text{ and } \quad 
v'(z)=\frac{\sqrt{y(v(z))(2\eps+y(v(z))}}{\eps+y(v(z))},
\end{equation}
while
the two boundary conditions in \eqref{notazioney} come from the fact that the monotone function $v$ has to attain the two equilibria $0$ and $1$ with zero derivative. 
In presence of a monostable reaction, 
the possible traveling speeds $c$ for the wave profiles solving \eqref{IIordMink} form an unbounded interval $[c^*_\eps, +\infty)$ (see Section \ref{sez2}), whose lower endpoint is known as \emph{critical speed}; the solution $v_\eps$ of \eqref{IIordMink} for $c=c_\eps^*$, unique up to translation in the variable $z=x+ct$, is called \emph{critical profile}. It is natural to wonder the behavior of $c_\eps^*$ and $v_\eps$ as $\eps \to 0^+$, that is, for vanishing diffusion. 

In the $0$-convection case, namely for constant $h$, it was shown in \cite{GarPP} that, under a mild additional assumption on $f$, the limit speed $\bar{c}=\lim_{\eps \to 0^+} c_\eps^*$ is strictly positive and coincides with $f(v_+)$, where $v_+ \in (0, 1)$ is the largest solution of the equation $\int_0^v f(s) \, ds=vf(v)$. Furthermore, given $\underline{v} \in (0, 1)$ and $\underline{z} \in \mathbb{R}$, if one denotes by  $\mathcal{V}_I^0(z; \underline{z}, \underline{v})$ the unique positive solution of 
$\displaystyle
\left\{
\begin{array}{l}
\bar{c}v'=f(v) \\ 
v(\underline{z})=\underline{v}
\end{array}
\right.
$ 
and defines 
$$
\mathcal{V}_L^0(z; \underline{z}, \underline{v})=\left\{
\begin{array}{ll}
0 & z \in (-\infty, \underline{z}-\underline{v}) \\ 
z-\underline{z}+\underline{v} & z \in [\underline{z}-\underline{v}, 1+\underline{z}-\underline{v}] \\ 
1 & z \in (1+\underline{z}-\underline{v}, +\infty), 
\end{array}
\right.
$$
it turns out that the limit profile for $\eps \to 0^+$ is obtained by gluing in a $C^1$-way such two functions, for a suitable choice of $\underline{z}$ and $\underline{v}$. Explicitly, for $\eps \to 0^+$ it holds that $v_\eps \to \bar{v}$ uniformly, where 
$$
\begin{array}{lll}
\bar{v}(z)
 = 
\left\{
\begin{array}{l} 
\mathcal{V}_L^0(z; 0, 1/2) \vspace{0.1cm} \\
\mathcal{V}_I^0(z; v_+-1/2, v_+) 
\end{array}
\right.
& 
\begin{array}{l} 
 z\in (-\infty, v_+-1/2] \\ 
 z \in (v_+-1/2, +\infty)
\end{array}
&
\text{if } v_+ \geq 1/2 \vspace{0.1cm}\\
\bar{v}(z)
= 
\left\{
\begin{array}{ll} 
\mathcal{V}_L^0(z; z_+, v_+)\vspace{0.1cm}\\ 
\mathcal{V}_I^0(z; 0, 1/2) 
\end{array}
\right.
& 
\begin{array}{l}
 z\in (-\infty, z_+] \\ 
z \in (z_+, +\infty)
\end{array}
& \text{if } v_+ < 1/2,
\end{array}
$$
being $z_+ < 0$ the unique real number for which $\mathcal{V}_I(z_+; 0, 1/2)=v_+$ \cite[Theorem 26]{GarPP}. Summarizing, the limit speed is \emph{strictly positive} and the limit profile \emph{becomes sharp on one side only}, near the equilibrium $0$.  This represents a deep difference with respect to models with a linear \cite{HilKim} or a saturating \cite{Gar21} diffusive term, for which $c_\eps^* \to 0$ for $\eps \to 0^+$ and the limit profile coincides with the Heaviside function. The crucial point leading to the result for the Born-Infeld operator consists in the possibility to provide a lower bound for the critical speed which is not compatible with a completely piecewise linear limit profile: roughly speaking, it was shown in \cite{GarPP} that $\bar{c} \geq \sup_{v \in (0, 1]} \dfrac{F(v)}{v} > F(1)$, while if it were $\bar{v}(z) =\mathcal{V}_L^0(z; 0, 1/2)$ for every $z \in \mathbb{R}$ one would have $\bar{c}=F(1)$ \cite[Proposition 23]{GarPP}.  

As is well known, the presence of a convective term may lead to deeply different outcomes with respect to the reaction-diffusion model: when the diffusive term is saturating, for example, it can provoke jumps in the traveling profiles \cite{GooKurRos} (see also \cite{GarStr}).
We can thus imagine that the conclusions drawn in \cite[Section 4]{GarPP} 
may drastically change for \eqref{IIordMink} and the solutions might possibly be led, for instance, to take the highest possible derivative for suitably chosen convections. 

In the case $f \equiv 0$, the presence of the sole convection steers in fact the limit profile for \eqref{IIordMink} to have derivative everywhere equal to $1$
when taking values different from the equilibria, as can be easily seen looking at the equivalent first-order two-point problem
\begin{equation}\label{soloconv}
\left\{
\begin{array}{l}
\displaystyle y'= (c+h'(v)) \frac{\sqrt{y(2\eps+y)}}{\eps+y} \vspace{0.1cm}\\
y(0)=0, \, y(1)= 0, \, y > 0 \text{ on } (0, 1).
\end{array}
\right.
\end{equation}
The differential equation herein is explicitly integrable; solving the associated Cauchy problem with initial condition $y(0)=0$ by separating variables yields the positive solution
$$
y_\eps(v)=\displaystyle \sqrt{\eps^2+ (cv+h(v))^2}-\eps,
$$
provided that $cv+h(v) > 0$ for every $v \in (0, 1)$.
Imposing the latter boundary condition in \eqref{soloconv}, one then finds that $c=-h(1)$, independently of $\eps$, and this implies that
$y_\eps(v) \to h(v) - h(1) v$ for every $v \in [0, 1]$ (actually, with uniform convergence). Due to \eqref{costruzionev}, this means that the limit profile $\bar{v}$ must necessarily have slope $1$ whenever it takes a value different from $0$ and $1$; since $v_\eps(0)=1/2$ for every $\eps > 0$, it will be $v_\eps'(z) \to 1$ for any $z \in (-1/2, 1/2)$ and $\bar{v}$ constant elsewhere, namely 
$\bar{v} \equiv \mathcal{V}_L(\cdot; 0, 1/2)$. 
Notice that, in order for \eqref{notazioney} to have a solution, here it has to be $h(v) > h(1) v$ for every $v \in (0, 1)$, that is, $h$ has to be concave. Convex convections support indeed the existence of \emph{decreasing} fronts, corresponding to problem \eqref{soloconv} where we take the opposite sign in the right-hand side. 
In any case, the purely convective model thus leads to a fully piecewise linear limit profile. 

It is then quite natural to wonder which scenarios can arise in presence of an interplay between reaction and convection and whether the above one-sided sharpening effect can be fully accentuated or destroyed in presence of sufficiently large concave or convex convections. 
Thanks to a careful study of problem \eqref{notazioney}, in Section \ref{sez3} we give a possible answer (Theorem \ref{thm-conv1}), after having provided an estimate of the critical speed associated with \eqref{IIordMink} (Section \ref{sez2}, to be compared, e.g., with \cite[Section 2]{GarPP}). The results are suitably complemented by the numerical simulations shown in Figures \ref{fig:1}--\ref{fig:4}. We point out that it would be interesting to investigate the shape of the limit critical profile in fourth-order models governed by the biharmonic operator, as well (see, e.g., \cite{BonSan}); in particular, one could wonder if 
the one-sided sharpening of the limit critical profile \cite{GarPP}, as well as the regularizing/sharpening effect of the convection in the present paper, could also arise in quasilinear versions of them.

\section{Estimating the critical speed}\label{sez2}

In this section, we focus on the bounds for the critical speed associated with \eqref{IIordMink}, dealing with its equivalent first-order reduction \eqref{notazioney}. 
We assume the following:
\begin{itemize}
\item[(F)] $f \in C([0, 1])$ is such that $f(0)=f(1)=0$, $f(s) > 0$ for $s \in (0, 1)$ and  
there exists $k > 0$ for which $f(s) \leq ks$, $ f(s) \leq k(1-s)$ for every $s \in [0, 1]$;
\item[(H)] $h \in C^2([0, 1])$ is such that $h(0)=h'(0)=0$. 
\end{itemize}
Assumption (F) is very common in the theory of traveling fronts and ensures their existence, as well as their proper behavior (namely, it excludes the possibility that they reach the equilibria in finite time). As for assumption (H), it only requires a standard behavior of $h$ at $0$, together with some basic regularity.  
\begin{definition}\label{def-adm}
We say that $c \in \mathbb{R}$ is an \emph{admissible speed} if problem \eqref{notazioney} has a solution. 
\end{definition}
Our first goal is to provide an estimate for the admissible speeds of increasing fronts. 
We first observe that, for such a kind of solutions, the condition 
$$
cv+h(v) > 0 \quad \text{ for every } v \in (0, 1]
$$
has to hold, as can be easily deduced by integrating equation \eqref{IIordMink} on $(-\infty, z)$, for any $z \in \mathbb{R} \cup \{+\infty\}$. Consequently, it has in particular to be $c+h(1) > 0$, in line with the discussion in the Introduction. Second, from $c > -h(v)/v$, passing to the limit for $v \to 0^+$ we infer 
$$
c \geq -h'(0) = 0,
$$
so that any admissible speed for \eqref{notazioney} is nonnegative.
Furthermore, from \eqref{notazioney} and the fact that $f$ is positive on $(0, 1)$ we deduce 
\begin{equation}\label{disugimportante}
y'(v) \leq (c+h'(v)-f(v)) \frac{\sqrt{y(v)(2\eps+y(v))}}{\eps+y(v)} \quad \text{ for every } v \in [0, 1], 
\end{equation}
since $\sqrt{y(2\eps+y)}\leq \eps + y$ for every $y \geq 0$; 
integrating \eqref{disugimportante} with the initial condition $y(0)=0$ necessarily implies $cv+h(v)-F(v) > 0$ 
for every $v \in (0, 1)$ and $y(v) \leq \sqrt{\eps^2+(cv+h(v)-F(v))^2}-\eps$, where we have set $F(v)=\int_0^v f(s) \, ds$.
Consequently, it has to be
\begin{equation}\label{estimate}
c \geq \sup_{v \in (0, 1]} \frac{F(v)-h(v)}{v}. 
\end{equation}
\smallbreak
The following lemma 
ensures the existence of an unbounded interval of admissible speeds.
\begin{lemma}\label{proponuova2}
Assume $(F)$ and $(H)$.
If 
\begin{equation}\label{Stimau}
c \geq \max_{v \in [0, 1]} f(v) - \min_{v \in [0, 1]} h'(v) + 2\sqrt{\eps\sup_{v \in (0, 1]}\frac{f(v)}{v}}, 
\end{equation}
then $c$ is admissible for \eqref{notazioney}. 
\end{lemma}
\begin{proof}
For fixed $\eps > 0$, we seek $\beta > 0$ such that the positive solution $y_\beta$ of 
$
\left\{
\begin{array}{l}
y_\beta'= \displaystyle \beta \frac{\sqrt{y_\beta(2\eps+y_\beta)}}{\eps+y_\beta} \vspace{0.1cm}\\
y_\beta(0)=0
\end{array}
\right.
$
satisfies
\begin{equation}\label{disug}
y_\beta'-(c+h'(v))\dfrac{\sqrt{y_\beta(2\eps+y_\beta)}}{\eps+y_\beta}+f(v) \leq 0 \quad \text{ for every } v \in [0, 1].
\end{equation}
This will provide a positive lower solution for the (forward) Cauchy problem associated with \eqref{notazioney}$_{1}$ with initial condition $y(0)=0$ (here \eqref{notazioney}$_{1}$ stands for the differential equation in \eqref{notazioney}). Since the unique solution of \eqref{notazioney}$_{1}$ fulfilling $y(1)=0$ cannot vanish in $(0, 1)$ due to the sign of $f$, this will imply the existence of a solution of \eqref{notazioney} by a standard uniqueness argument. 
Computing now the explicit expression of $y_\beta$, \eqref{disug}
will be true if, for every $v \in (0, 1]$, it holds
$$
\beta^2 -(c+h'(v)-f(v))\beta + \eps \frac{f(v)}{v}\leq 0.
$$
This inequality is implied by the condition 
$$
\beta^2 -\Big(c+\min_{v \in [0, 1]} h'(v)-\max_{v \in [0, 1]} f(v)\Big)\beta + \sup_{v \in (0, 1]} \eps\frac{f(v)}{v} \leq 0;
$$
thanks to \eqref{Stimau}, there exists at least one positive $\beta$ for which such an inequality holds. One can choose, for instance, 
$$
\beta=\frac{c+\min_{v \in [0, 1]} h'(v)-\max_{v \in [0, 1]} f(v)}{2} > 0,
$$
correspondingly finding the desired positive lower solution $y_\beta$.
\end{proof}
We point out that the value appearing in the right-hand side of \eqref{Stimau} is strictly positive, coherently with the previous discussion, since $\min_{v \in [0, 1]} h'(v) \leq h'(0) =0$. Now, 
the argument in the proof of Lemma \ref{proponuova2} ensures that if $c$ is admissible and $c' > c$, then also $c'$ is admissible; moreover, arguing as in the proof of \cite[Proposition 3.2]{CoeSan}, the use of Arzel\`a-Ascoli Theorem ensures that also the infimum of admissible speeds is an admissible speed. The set of admissible speeds is then an interval, denoted by $[c_\eps^*, +\infty)$; its lower endpoint $c_\eps^*$ is called \emph{critical speed} and the associated front profile is called \emph{critical profile}. Notice that Lemma \ref{proponuova2} provides an upper bound for $c_\eps^*$. Moreover, we underline that $c_\eps^*$ is monotone increasing in $\eps$, as a result of comparison principles applied to the first-order backward problem for \eqref{notazioney}$_{1}$ with initial condition $y(1)=0$, which enjoys uniqueness.

We are interested in the limit behavior of $c_\eps^*$ and $v_\eps$ for $\eps \to 0^+$;
to this end, we will examine the solution 
$y_\eps$ of the first-order reduction \eqref{notazioney}.
We set
$$
\bar{c}=\lim_{\eps \to 0^+} c_\eps^*, \quad \bar{v}=\lim_{\eps \to 0^+} v_\eps, \quad \bar{y}=\lim_{\eps \to 0^+} y_\eps.
$$
Notice indeed that $c_\eps^*$ converges since it is monotone in $\eps$. 
On the other hand, the front profiles $v_\eps$ are equi-Lipschitz continuous and bounded in $C^1$ since $0 \leq v_\eps \leq 1$ and $\Vert v_\eps' \Vert_{L^\infty(\mathbb{R})} \leq 1$, hence there exists a Lipschitz continuous function $\bar{v}$ such that  $\bar{v}(z)=\lim_{\eps \to 0^+} v_\eps(z)$ for every $z \in \mathbb{R}$ and $v_\eps \to \bar{v}$ uniformly on compact subsets of $\mathbb{R}$ (actually, the convergence is uniform on the whole real line in view of \cite[Lemma 2.4]{Diek}). Finally, for what concerns $y_\eps$, due to the sign of $c_\eps^*$ for sure it will be $y_\eps' \geq -\max_{v \in [0, 1]} \vert h'(v) \vert - f(v)$, hence the backward solution of $y'=-\max_{v \in [0, 1]} \vert h'(v) \vert-f(v)$ satisfying $y(1)=0$ is a positive upper solution for $y_\eps$, for every $\eps$. Consequently, 
there exists $Y_\infty > 0$ such that 
\begin{equation}\label{y1}
\Vert y_\eps \Vert_{L^\infty(0, 1)} \leq Y_\infty
\end{equation} 
for every $\eps > 0$; 
since $y \mapsto \dfrac{\sqrt{y(2\eps+y)}}{\eps + y}$ is increasing, the differential equation in \eqref{notazioney} then yields the existence of $Y_\infty' >0$ for which  
\begin{equation}\label{y2}
\Vert y_\eps' \Vert_{L^\infty(0, 1)} \leq Y_\infty'
\end{equation}
for every $\eps > 0$. 
By the Arzelà-Ascoli Theorem, the bounds \eqref{y1} and \eqref{y2} imply that $y_\eps$ converges uniformly to its limit $\bar{y}$ for $\eps \to 0^+$ and $\bar{y}$ is a continuous function, nonnegative on $[0, 1]$. 

To begin with, we give an estimate of the limit critical speed $\bar{c}$ for $\eps \to 0^+$ in case the limit profile contains a linear piece with slope $1$.  
\begin{proposition}
Assume that, for some real numbers $z_0 < z_1$, $v_0 \in [0, 1)$, it holds $\bar{v}(z)=\mathcal{V}_L(z; z_0, v_0)$ on $[z_0, z_1]$
and let $\bar{v}(z_1)=v_1\in (0, 1]$ (so that $\bar{v}(z)=z-z_0+v_0$ on $[z_0, z_1]$, being $z_1-z_0=v_1-v_0$). Then, 
$$
\bar{c}=\frac{F(v_1)-F(v_0)-h(v_1)+h(v_0)}{v_1-v_0}+\frac{\bar{y}(v_1)-\bar{y}(v_0)}{v_1-v_0}.
$$
\end{proposition}
\begin{proof}
Since $\{v_\eps\}_\eps$ is bounded in $C^1([z_0, z_1])$, the family $\{v_\eps'\}_\eps$ is bounded in $L^2(z_0, z_1)$ and hence there exists $w \in L^2(z_0, z_1)$ for which, up to subsequences, $v_\eps' \rightharpoonup w$ in $L^2(z_0, z_1)$ for $\eps \to 0^+$. Passing to a further subsequence, if necessary, $w$ has to coincide almost everywhere with $\bar{v}'$ in the points where $\bar{v}$ is differentiable, so that $w(z)=1$ for almost every $z \in [z_0, z_1]$. Multiplying the differential equation in \eqref{IIordMink} by $v'_\eps$ and integrating on $[z_0, z_1]$ one obtains, for every $\eps > 0$,   that
\begin{eqnarray*}
c_\eps^* \Vert v_\eps' \Vert^2_{L^2(z_0, z_1)}\!+ \!\!\int_{z_0}^{z_1} h'(v_\eps(z)) v_\eps'(z)^2 \, dz 
 & \!\!\!\!= & \!\!\!\!  \eps\left(\frac{1}{\sqrt{1-v_\eps'(z_1)^2}}\!-\!\frac{1}{\sqrt{1-v_\eps'(z_0)^2}}\right) \!\!+\!F(v_\eps(z_1))\!-\!F(v_\eps(z_0)) \nonumber \\
&\!\!\! = & \!\!\!\! y_\eps(v_\eps^1)-y_\eps(v_\eps^0) + F(v_\eps^1)-F(v_\eps^0), 
\end{eqnarray*}
where $z_\eps$ is the inverse function of $z \mapsto v_\eps(z)$ and $v_\eps^0:=v_\eps(z_0)$, $v_\eps^1:=v_\eps(z_1)$. 
Using the weak semicontinuity of the norm, the first summand in the left-hand side is bounded from below by $\bar{c}(v_1-v_0)$ in the limit for $\eps \to 0^+$. 
As for the second one, we observe that $h'(v_\eps)v_\eps' \to h'(\bar{v})w$ strongly in $L^2(z_0, z_1)$, in view of the uniform convergence of $v_\eps$ to $\bar{v}$ and of the fact that $h \in C^2$. Since $v_\eps' \rightharpoonup w$ in $L^2(z_0, z_1)$, we conclude that
$$
\int_{z_0}^{z_1} h'(v_\eps(z)) v_\eps'(z)^2 \, dz  \to \int_{z_0}^{z_1} h'(\bar{v}(z)) w(z)^2 \, dz = h(v_1)-h(v_0), 
$$
where the last equality follows from the properties of absolutely continuous functions.
Being $\lim_{\eps \to 0^+} v_\eps^0 = v_0$ and $\lim_{\eps \to 0^+ }v_\eps^1 = v_1$, using \eqref{y1} and \eqref{y2} one has that $y_\eps(v_\eps^0) \to \bar{y}(v_0)$ and $y_\eps(v_\eps^1) \to \bar{y}(v_1)$ for $\eps \to 0^+$. Therefore, we infer 
$$
\bar{c} \leq \frac{\bar{y}(v_1)-\bar{y}(v_0)}{v_1-v_0} + \frac{F(v_1)-F(v_0)-h(v_1)+h(v_0)}{v_1-v_0}. 
$$
On the other hand, integrating the inequality \eqref{disugimportante} - with $c=c_\eps^*$ -  between 
$v_0$ and $v_1$ provides $\displaystyle \sqrt{y_\eps(v_1)(2\eps+y_\eps(v_1))}-\sqrt{y_\eps(v_0)(2\eps+y_\eps(v_0))} \leq c_\eps^*(v_1-v_0)+h(v_1)-h(v_0)-F(v_1)+F(v_0)$, which passing to the limit for $\eps \to 0^+$ yields the reversed inequality for $\bar{c}$.
The conclusion follows.
\end{proof}
\begin{corollary}\label{conv1parz}
Let $\bar{v}(z) = \mathcal{V}_L(z; 0, 1/2)$ for every $z \in \mathbb{R}$. 
Then, $\bar{c}=F(1)-h(1)$. 
\end{corollary}

\section{The role of the convective term in the shape of the limit profile }\label{sez3}

From the discussion in Section \ref{sez2}, 
we have understood that a key quantity in determining the asymptotic value of the critical speed and the asymptotic shape of the critical profile is 
$$
S(v):=\dfrac{F(v)-h(v)}{v}.
$$
In view of assumptions (F) and (H), $S$ can be extended by continuity setting $S(0)=0$, and for this reason $\sup_{v \in (0, 1]} S(v)  \geq 0$. 
In the next statement, we see how different assumptions on $S$ produce different outcomes regarding $\bar{c}$ and $\bar{v}$; the exposition is maintained to a more readable level in order to avoid overloading the contents and the proofs, since we are mainly interested in highlighting the different phenomena which can arise rather than in providing the corresponding optimal assumptions. In any case, our results fully cover the most natural case of Fisher-Burgers type equations (see \eqref{esempio}), for which we are able to give a complete picture.  \\ In the following, we set $\mathcal{V}_L(z; \underline{z}, \underline{v})=\mathcal{V}_L^0(z; \underline{z}, \underline{v})$, and
we denote by $\mathcal{V}_I(\cdot; \underline{z}, \underline{v})$ 
the unique $C^1$-solution of 
\begin{equation}\label{pinviscid}
\left\{
\begin{array}{l}
(\bar{c}+h'(v))v'=f(v) \\ 
v(\underline{z})=\underline{v};
\end{array}
\right.
\end{equation}
here we implicitly assume that $\bar{c}+h'(\underline{v}) > 0$, in line with the monotonicity of the solutions we are interested in and in order to avoid degeneracy. Of course, the solution of \eqref{pinviscid} might possibly be defined only on a neighborhood of $\underline{z}$ or might escape the interval $[0, 1]$, depending on the expressions of $f$ and $h'$. We still comment about this after Theorem \ref{thm-conv1}. 

With these preliminaries, we now have the following. 
\begin{theorem}\label{thm-conv1}
Let $f$ and $h$ fulfill assumptions (F) and (H). The following hold:
\begin{itemize}
\item[1)] if $S'(v) < 0$ for every $v \in (0, 1]$, then $\bar{c}=0$ and $\bar{v}\equiv \mathcal{V}_I(\cdot; 0, 1/2)$ as long as $\bar{v} > 0$, while $\bar{v}=0$ elsewhere;
\item[2)] if $0 < \sup_{v \in (0, 1]} S(v) \neq S(1)$, then $\bar{c}=f(v_+)-h'(v_+)$, where $v_+$ is the largest solution of the equation $F(v)-h(v)=v(f(v)-h'(v))$; moreover, if $f-h'$ has a unique maximum point in $[0, 1]$, the limit profile $\bar{v}$ is given by  
$$
\begin{array}{lll}
\bar{v}(z)
 = 
\left\{
\begin{array}{l} 
\mathcal{V}_L(z; 0, 1/2)  \\
\mathcal{V}_I(z; v_+-1/2, v_+) 
\end{array}
\right.
& 
\begin{array}{l} 
 z\in (-\infty, v_+-1/2] \\ 
 z \in (v_+-1/2, +\infty)
\end{array}
&
\text{if } v_+ \geq 1/2 \vspace{0.1cm}\\
\bar{v}(z)
= 
\left\{
\begin{array}{ll} 
\mathcal{V}_L(z; z_+, v_+)\\ 
\mathcal{V}_I(z; 0, 1/2) 
\end{array}
\right.
& 
\begin{array}{l}
 z\in (-\infty, z_+] \\ 
z \in (z_+, +\infty)
\end{array}
& \text{if } v_+ < 1/2,
\end{array}
$$
where, in the latter case, $z_+ < 0$ is the unique real number for which $\mathcal{V}_I(z_+; 0, 1/2)=v_+$;
\item[3)] if $S'(v) > 0$ for every $v \in (0, 1)$, then $\bar{c}=F(1)-h(1)$ and $\bar{v}\equiv \mathcal{V}_L(\cdot; 0, 1/2)$. 
\end{itemize}
\end{theorem}
\begin{proof}
We first notice that if $\bar{c} > 0$, then there exists $v_1 \in (0, 1)$ such that $\bar{y}(v) > 0$ for every $v \in (0, v_1]$; indeed, if $c < \bar{c}$ is sufficiently small, then the positive solution of \eqref{notazioney}$_{1}$ satisfying $y(0)=0$ is greater or equal than the positive solution of 
$$
\left\{
\begin{array}{l}
\displaystyle y'= c \frac{\sqrt{y(2\eps+y)}}{\eps+y} \vspace{0.1cm}\\
y(0)=0,
\end{array}
\right.
$$
uniformly in $\eps$.
This can be seen, for instance, performing an argument similar to the one in the proof of \cite[Lemma 27]{GarPP}, thanks to the fact that $h'(0)=0$ and hence $\bar{c}+h'(v) > 0$ in a right neighborhood of $0$. Hence, it makes sense to define $\widetilde{v}:=\sup\{v \in (0, 1] \mid \bar{y} > 0 \text{ on } (0, v)\}$. One has that $\bar{y}(\widetilde{v})=0$, by the continuity of $\bar{y}$; moreover, for every $[\alpha, \beta] \subset (0, v_1)$, it holds that $y_\eps \to \bar{y}$ in $C^1([\alpha, \beta])$, where $\bar{y}(v)=\bar{c}v+h(v)-F(v)$ (as can be seen letting $\alpha \to 0^+$). 

The proof then proceeds differently according to the considered case.
\begin{itemize}
\item[1)] If by contradiction it were $\bar{c} > 0$, by the previous observations one would have $\bar{y}(\widetilde{v})=0$ and hence, by the expression of $\bar{y}$, 
$$
\bar{c}=\frac{F(\widetilde{v})-h(\widetilde{v})}{\widetilde{v}} \leq 0, 
$$
which is a contradiction. Hence, $\bar{c}=0$. For fixed $\xi \in (0, 1)$, let now $M_\eps=\max_{v \in [\xi, 1]} y_\eps(v)$ and let $v_{M, \eps} \in [\xi, 1)$ be such that $y_\eps(v_{M,\eps})=M_\eps$ (we drop the dependences on $\xi$ for the sake of readability); up to subsequences, one can assume that $v_{M, \eps} \to v^* \in [\xi, 1]$. If by contradiction it were $\eps/M_\eps \to 0$, passing to the limit in 
$$
0 \geq y_\eps'(v_{M, \eps})=(c_\eps^*+h'(v_{M, \eps})) \frac{\sqrt{M_\eps(2\eps+M_\eps)}}{\eps+M_\eps}-f(v_{M, \eps})
$$
one would obtain $h'(v^*) \leq f(v^*)$, implying that $S'(v^*) \geq 0$ against the assumption. It follows that $y_\eps \to 0$ uniformly and, for any fixed $\xi > 0$, $M_\eps \to 0$ with order $\eps$ (if this occurred with a stronger order, one would easily find the contradiction $y_\eps' \to -f < 0$). For fixed $\zeta \in \mathbb{R}$, using the expression of $y$ provided by \eqref{costruzionev} then ensures the existence of $0 < K  < 1$ such that $\Vert v'_\eps \Vert_{L^\infty(\zeta, +\infty)} \leq K$.
For any $\psi \in C_c^\infty([\zeta, +\infty))$, one can then pass to the limit for $\eps \to 0^+$ in
$$
-\int_{\zeta}^{+\infty} \!\!\frac{\eps v_{\eps}'(z)}{\sqrt{1-(v_{\eps}'(z))^2}} \psi'(z) \, dz + \int_{\zeta}^{+\infty}\!\! (c_{\eps}^*v_\eps(z)+h(v_\eps(z))) \psi'(z) \, dz + \int_{\zeta}^{+\infty} \!\! f(v_{\eps}(z)) \psi(z) \, dz = 0,
$$
finally obtaining that $\bar{v}$ satisfies $h'(\bar{v})\bar{v}'-f(\bar{v})=0$ whenever $\bar{v} > 0$ (that is, for $\zeta > -\infty$). It follows that $\bar{v}(z) = \mathcal{V}_I(z; 0, 1/2)$ for every $z \in \mathbb{R}$ such that $\bar{v}(z) > 0$.
\item[2)] In this second case, we have $\bar{c} > 0$ in view of \eqref{estimate} and thus $\widetilde{v} > 0$ is well defined. 
Moreover, we observe that necessarily $\widetilde{v} < 1$, otherwise $\bar{c}=F(1)-h(1)$ and hence $\bar{y}(v)=(F(1)-h(1))v-F(v)+h(v) > 0$ for every $v \in (0, 1)$, contradicting the assumption $\sup_{v \in (0, 1]} S(v) \neq S(1)$. 
We now claim that 
\begin{equation}\label{valuebarc}
\bar{c}=\frac{F(\widetilde{v})-h(\widetilde{v})}{\widetilde{v}}=f(v_m)-h'(v_m)=f(\widetilde{v})-h'(\widetilde{v}), 
\end{equation}
where $v_m \in (0, \widetilde{v})$ is such that $\bar{y}(v_m)=M:=\max_{v \in [0, \widetilde{v}]} \bar{y}(v)$. 
The first equality in \eqref{valuebarc} is a straight consequence of the fact that $\bar{y}(\widetilde{v})=0$, 
while the second one follows from the $C^1$-convergence of $y_\eps$ to $\bar{y}$ on any interval $[\alpha, \beta] \subset (0, \widetilde{v})$, since
$$
0 \leftarrow \displaystyle y_\eps'(v_m)= (c_\eps^*+h'(v_m)) \frac{\sqrt{y_\eps(v_m)(2\eps+y_\eps(v_m))}}{\eps+y_\eps(v_m)} - f(v_m) \to \bar{c}+h'(v_m) - f(v_m),
$$
being $y_\eps(v_m) \to M > 0$. As for the last equality in \eqref{valuebarc}, we first notice that, for fixed $v_1 > \widetilde{v}$, integrating the inequality $y_\eps'(v) \leq (c_\eps^*+h'(v)-f(v)) \dfrac{\sqrt{y_\eps(v) (2\eps+y_\eps(v))}}{\eps+y_\eps(v)}$ between $\widetilde{v}$ and $v_1$ yields 
$
\sqrt{y_\eps(v_1)(2\eps+y_\eps(v_1))}-\sqrt{y_\eps(\widetilde{v})(2\eps+y_\eps(\widetilde{v}))} \leq c_\eps^*(v_1-\widetilde{v})+h(v_1)-h(\widetilde{v}) - F(v_1)+F(\widetilde{v}),
$ 
which passing to the limit for $\eps \to 0^+$ produces (recalling that $\bar{y}(\widetilde{v})=0$ and $\bar{y}(v_1) \geq 0$)
$$
\bar{c} \geq \frac{F(v_1)-F(\widetilde{v})-h(v_1)+h(\widetilde{v})}{v_1-\widetilde{v}};
$$
passing now to the limit for $v_1 \to \widetilde{v}^+$, we obtain $\bar{c} \geq f(\widetilde{v})-h'(\widetilde{v})$.
As for the reversed inequality, we observe that since $\bar{y}(\widetilde{v}) = 0 < \bar{y}(v)$ for every $v \in (0, \widetilde{v})$, there exists a sequence $v_n \nearrow \widetilde{v}$ such that $\bar{y}'(v_n) \leq 0$. 
Recalling the expression of $\bar{y}$ in the interval $(0, \widetilde{v})$ and the $C^1$ convergence of $y_\eps$ to $\bar{y}$, it then follows that
$$
0 \geq \bar{y}'(v_n) = \bar{c}+h'(v_n) - f(v_n), 
$$
which yields the desired inequality 
passing to the limit for $n \to +\infty$. 
We now proceed similarly as in case 1) and as in the proof of \cite[Theorem 26]{GarPP}. For fixed $\sigma > 0$, we define $M_\eps=\max_{v \in [\widetilde{v}+\sigma, 1]} y_\eps(v)$ and we let $v_{M, \eps} \in [\widetilde{v}+\sigma, 1)$ be such that $y_\eps(v_{M, \eps})=M_\eps$.  
Up to subsequences, $v_{M, \eps} \to v^* \in [\widetilde{v}+\sigma, 1]$; if it were $\eps/M_\eps \to 0$, arguing as in case 1) one would find $\bar{c}+h'(v^*)-f(v^*) \leq 0$, which is impossible in view of the assumption on $f-h'$ (since the equation $\bar{c}+h'(v)-f(v)=0$ has the two distinct solutions $v_m$ and $\widetilde{v}$). Hence $M_\eps \to 0$ with order $\eps$, $y_\eps \to 0$ uniformly in $[\widetilde{v}+\sigma, 1]$ for every $\sigma > 0$, $\widetilde{v}=v_+$ and $\bar{v}$ is piecewise linear with slope $1$ as long as it takes values in $(0, v_+)$ and coincides with the solution of \eqref{pinviscid} gluing in a $C^1$-way while taking values in $(v_+, 1)$. 
This is equivalent to the statement in case 2).
\item[3)] Since $h'(0)=0$, the assumption $S' > 0$ and \eqref{estimate} imply that 
$\bar{c} \geq \sup_{v \in (0, 1]} S(v)=F(1)-h(1) > 0$. 
With the same notation as above, we here have that $\widetilde{v}$ necessarily coincides with $1$, for if by contradiction $\widetilde{v} < 1$, one would have 
$\bar{y}(\widetilde{v})=\bar{c}\widetilde{v}+h(\widetilde{v})-F(\widetilde{v})=0$, implying, in view of the fact that $S' > 0$, the contradiction $\bar{c}= \dfrac{F(\widetilde{v})-h(\widetilde{v})}{\widetilde{v}} < F(1)-h(1)$.
Hence, $\bar{c}=F(1)-h(1)$ and $\bar{y}(v)=\bar{c}v+h(v)-F(v) > 0$ for every $v \in (0, 1)$. Using a similar argument as the one mentioned in the Introduction for the $0$-convection case, this implies that $\bar{v} \equiv \mathcal{V}_L(\cdot; 0, 1/2)$, concluding the statement.  
\end{itemize}
\end{proof}
Some comments are in order, in particular about the above case 1). Setting $z_0=\inf\{z \in (-\infty, 0) \mid \mathcal{V}_I(z; 0, 1/2) > 0\}$, it may be $z_0=-\infty$ or $z_0 \in (-\infty, 0)$, according to whether the improper integral $\int_0^{1/2} h'(s)/f(s) \, ds$ diverges or not. In the former case, $\mathcal{V}_I(\cdot; 0, 1/2)$ (henceforth briefly denoted by $\mathcal{V}_I$) is positive and regular on the whole real line and $v_\eps \to \bar{v}$ with $C^1$-convergence, while in the latter one it holds $\mathcal{V}_I(z_0)=0$ and the differential equation for $\mathcal{V}_I$ degenerates at $z_0$. If $\mathcal{V}_I'(z_0) > 0$, $\mathcal{V}_I$ may then be prolonged so as to escape the interval $[0, 1]$, while if $\mathcal{V}_I'(z_0)=0$ it may be prolonged in a $C^1$ way to $0$ in the interval $(-\infty, z_0)$. Noticing that $v_\eps \to \mathcal{V}_I$ in $C^1(z, +\infty)$ for every $z > z_0$, similarly as in \cite[Remark 30]{GarPP}, it will be
$$
\lim_{z \to z_0^+} \bar{v}'(z)=\lim_{z \to z_0^+} \mathcal{V}_I'(z)=\lim_{z \to z_0^+} \frac{f(\mathcal{V}_I(z))}{h'(\mathcal{V}_I(z))};
$$
hence, a crucial role is here played by the limit $\ell:=\lim_{s \to 0^+} f(s)/h'(s)$. If $\ell=0$, then $\bar{v}$ may be prolonged to a globally $C^1$-function, being $\lim_{z \to z_0^+}\bar{v}'(z)=0=\lim_{z \to z_0^-} \bar{v}'(z)$ (this last equality holding since necessarily $\bar{v} \equiv 0$ in $(-\infty, z_0)$). Otherwise, if $\ell > 0$ then $\bar{v}$ is sharp near $0$. 
This issue does not arise for the inviscid piece of limit profile in case 2), since the boundedness of the considered front profiles in $C^1$, together with the fact that $f$ is positive on $(0, 1)$, prevents $c+h'(v)$ from vanishing, otherwise one would reach the contradiction $0=f(v)$ for some $v \in [v_+, 1)$. 

We make some further remarks about the statement of Theorem \ref{thm-conv1}. 
We notice that case 1) only occurs if $h'(s) > 0$ for every $s \in (0, 1]$, namely only for suitable \emph{strictly increasing convections}. Indeed, the conditions $S(0)=0$ and $S'(v) < 0$ for every $v \in (0, 1]$ imply that $S(v) < 0$ for every $v> 0$; since moreover 
$$
S'(v)=\frac{f(v)-h'(v)-S(v)}{v}, 
$$
it has necessarily to be $f(v)-h'(v)-S(v) < 0$, that is, $f(v)+\vert S(v) \vert < h'(v)$ for $v > 0$. Notice that this implies that the differential equation in \eqref{pinviscid} can degenerate only when $v=0$. 

Case 3) \emph{cannot} instead occur if $h$ is \emph{everywhere increasing}, since from the fact that $S'(v) > 0$ for every $v \in (0, 1)$ it follows that  $f(1)-h'(1)-S(1)=-h'(1)-S(1) \geq 0$, hence $h'(1)$ has to be strictly negative (in this case $S(1) > 0$). However, the assumption $S'(v) > 0$ for every $v \in (0, 1)$ may be fulfilled for convections which are \emph{locally} increasing near $0$, like $h(s)=s^2(\delta-s)$, for a sufficiently small $\delta >0$, provided that $f(v) >h'(v)+S(v)$ for every $v \in (0, 1)$. 

In particular, since $h'(0)=0$, \emph{concave} and \emph{convex} convections are ruled out, respectively, from cases 1) and 3). 
Case 2) depends instead on the quantitative interplay between $f$ and $h'$, rather than on the sign of $h'$, hence in principle it can occur regardless of the monotonicity (and of the convexity) of $h$. One may wonder if the sign assumption $S' < 0$ (resp., $S' > 0$) on the derivative of $S$  could be replaced by a weaker condition like $\sup_{v \in (0, 1]} S(v)=0$ (resp., $\sup_{v \in (0, 1]} S(v)= F(1)-h(1)$), but for the sake of brevity we have preferred a slightly stronger hypothesis, in order to proceed with a simpler and shorter proof. 

\begin{figure}[h!]
\centering
\includegraphics[scale=0.46]{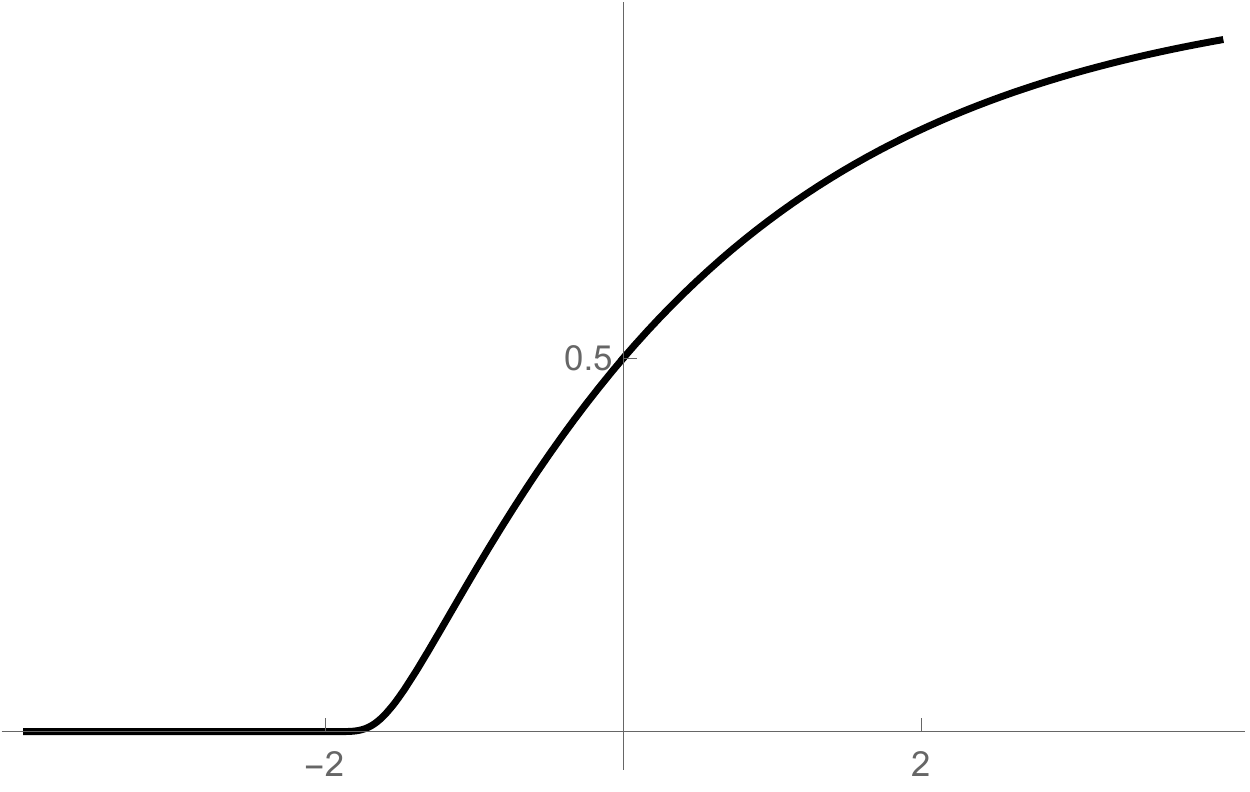}
    \hspace{1cm} \includegraphics[scale=0.49]{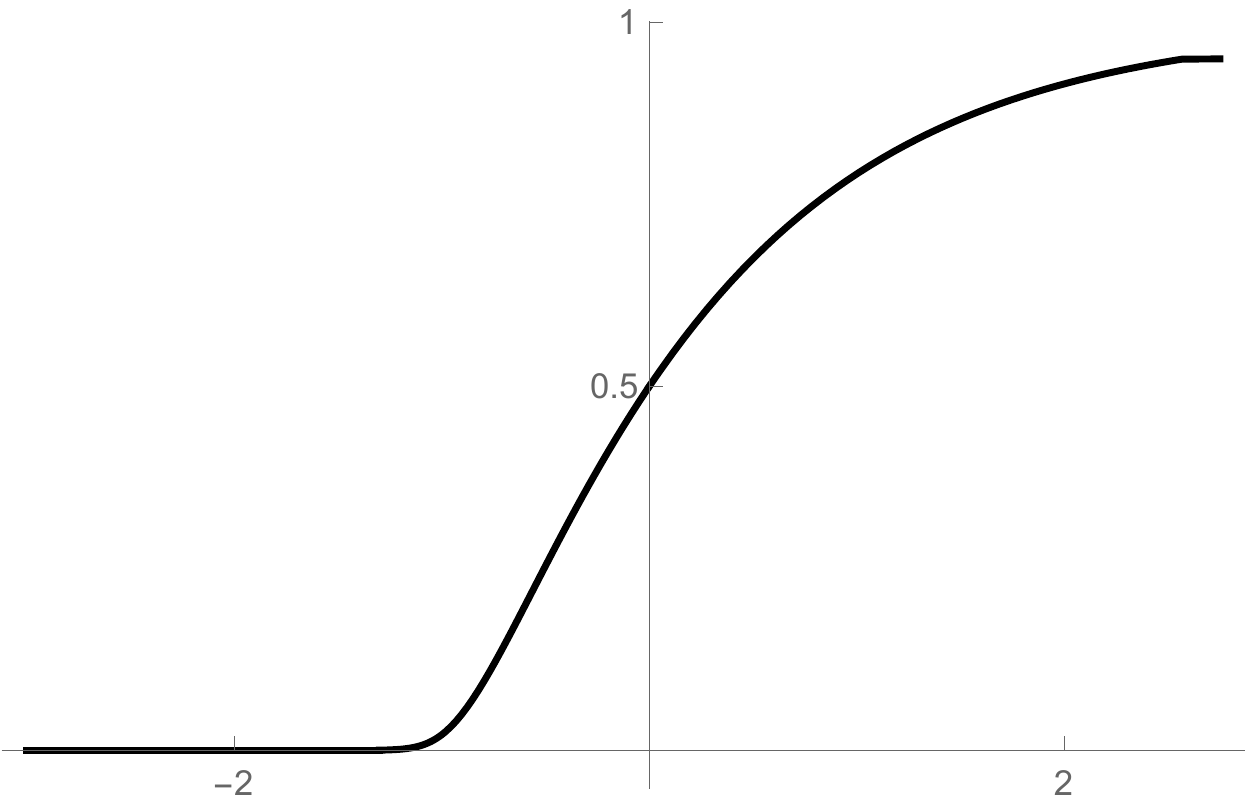}\\
    \includegraphics[scale=0.46]{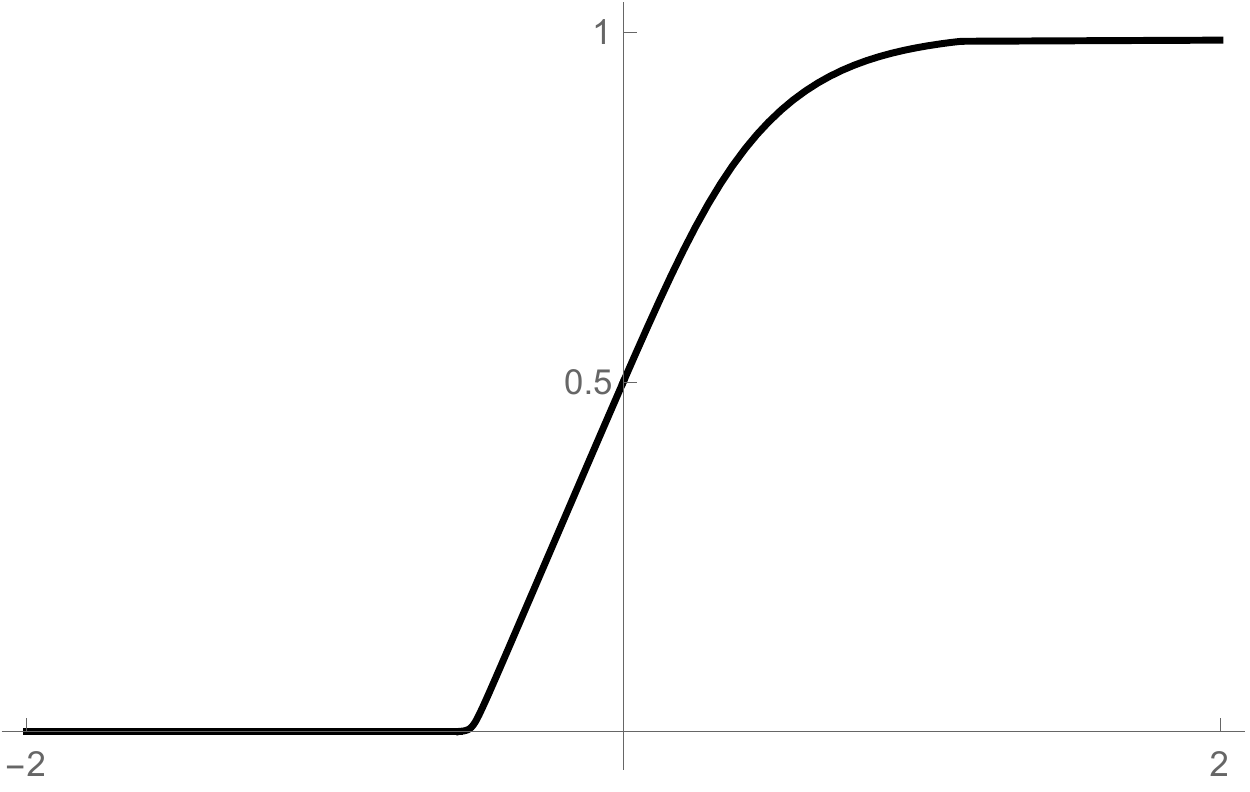}
\hspace{1cm} \includegraphics[scale=0.49]{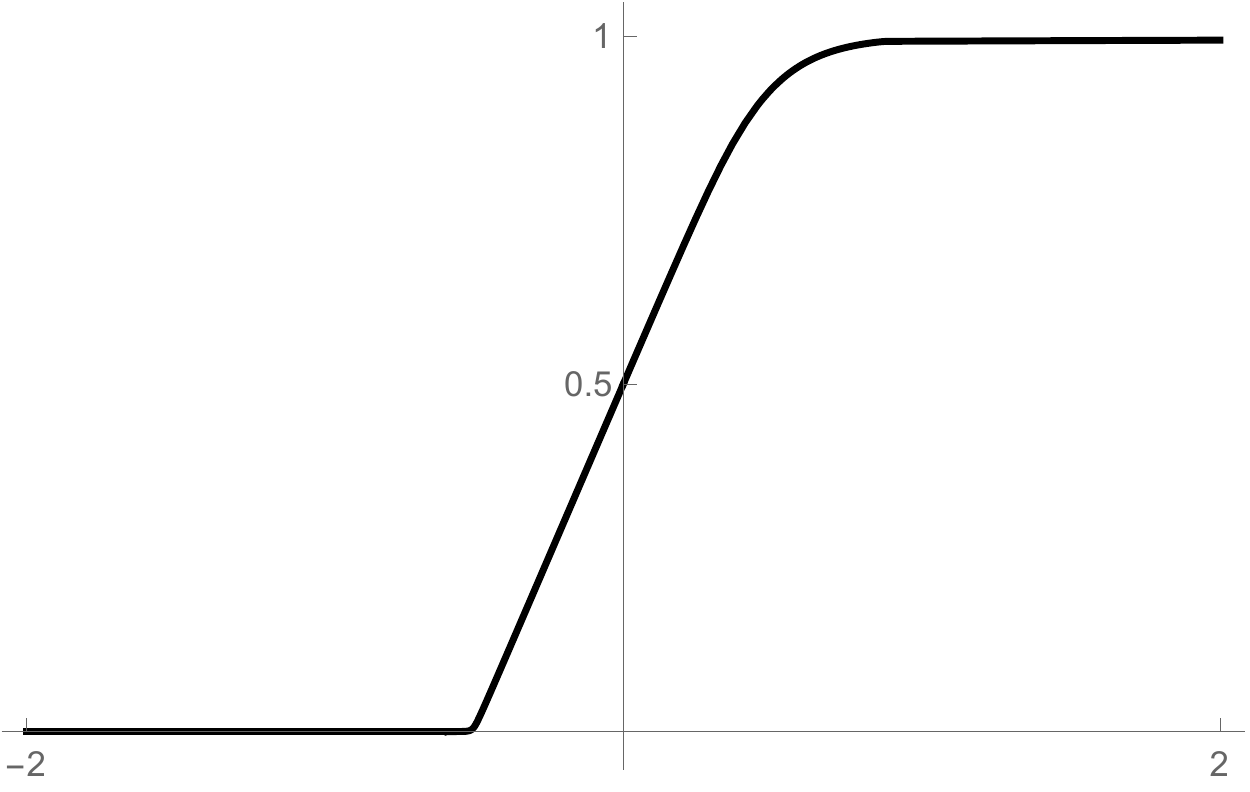}\\
     \includegraphics[scale=0.46]{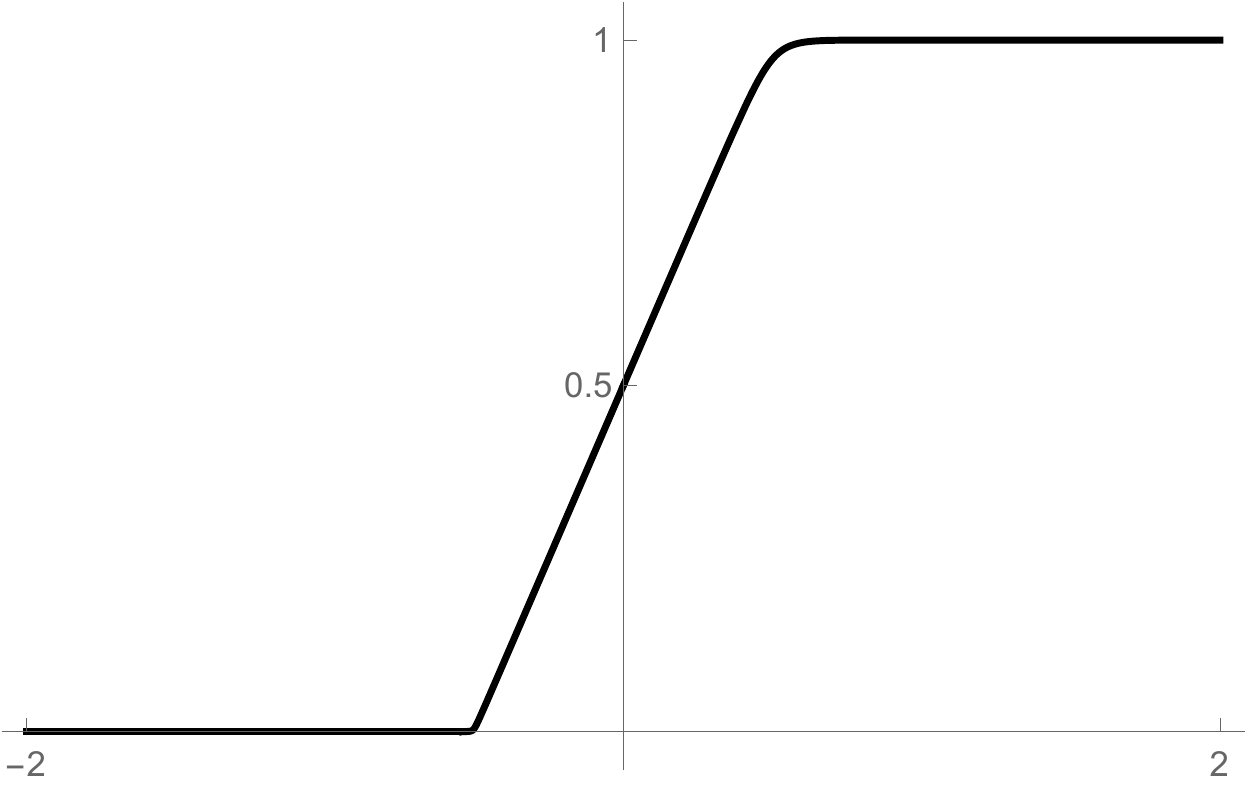}
    \hspace{1cm} \includegraphics[scale=0.49]{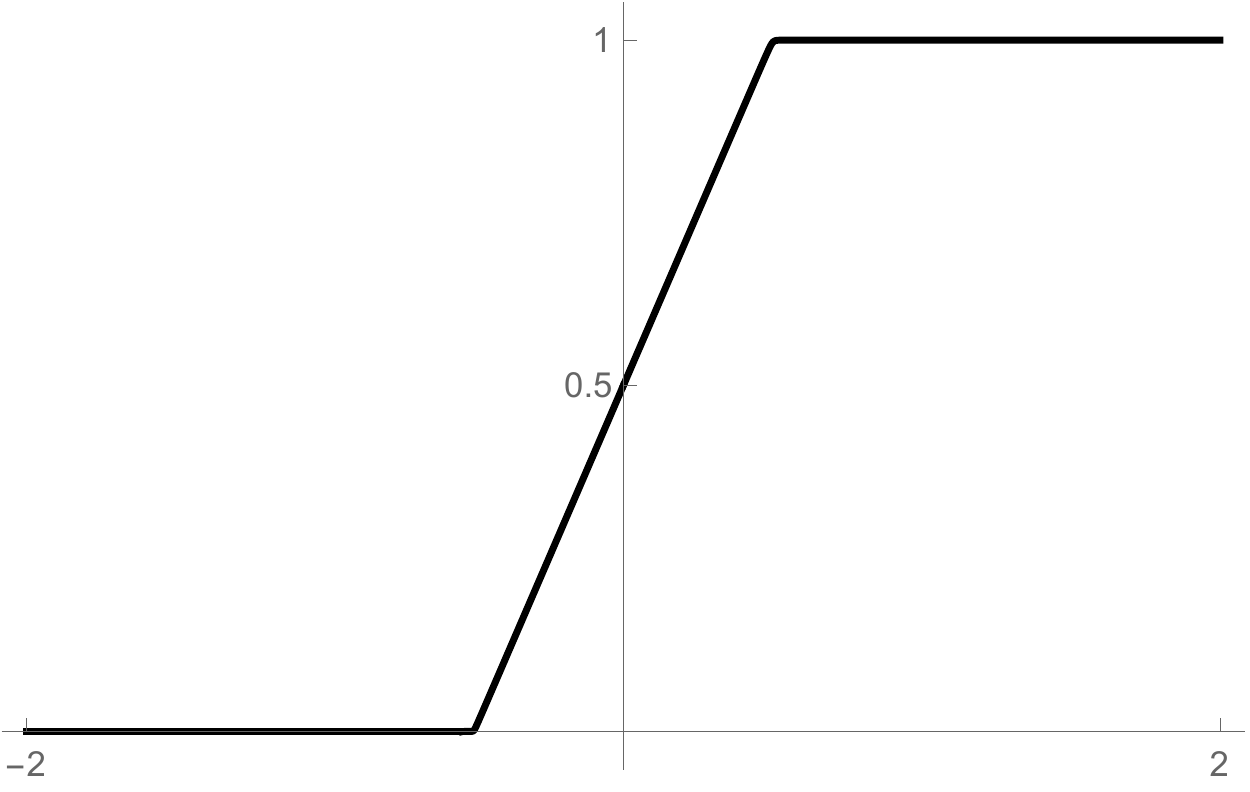}\\
\caption{For $\eps, k, \alpha$ as in \eqref{sceltefissate}, we show the graph of the critical profile $v_\eps$ for \eqref{esempio}. The values of $\bar{c}$ predicted by Theorem \ref{thm-conv1} are, respectively, $\bar{c}=0$ (top line), $\bar{c}=0.152, 0.187$ (middle line), $\bar{c}=1/3, 2/3$ (bottom line), while here we have found the approximations $c_\eps^* \approx 0.07, 0.09$ (top line), $c_\eps^* \approx 0.163, 0.234$ (middle line), $c_\eps^* \approx 0.336, 0.667$ (bottom line).\vspace{2cm}}\label{fig:1}
\end{figure}

We illustrate the statement of Theorem \ref{thm-conv1} for the Fisher-Burgers type equation 
\begin{equation}\label{esempio}
u_t=\eps \left(\frac{u_x}{\sqrt{1-u_x^2}}\right)_x  - (\alpha u^2)_x + ku(1-u) =0,
\end{equation}
in dependence on two parameters $\alpha\in \mathbb{R}$, $k > 0$, showing in Figure \ref{fig:1} some numerical simulations obtained through the numerical integration of \eqref{notazioney}, using Wolfram Mathematica$^{\copyright}$ software. To this end, we set
\begin{equation}\label{sceltefissate}
\begin{array}{ll}
\eps=2\cdot 10^{-3}, k=1 \text{ and: } \alpha=1, \; \alpha=0.5(=k/2) \text{\,(top line)}, \vspace{0.1cm}\\
\alpha=0.05, \; \alpha=-0.05 \text{\,(middle line)}, \, \alpha=-1/6(=-k/6), \; \alpha= -0.5 \text{\,(bottom line)}. 
\end{array}
\end{equation}
For equation \eqref{esempio} one has $S(v)=\displaystyle \left(\frac{k}{2}-\alpha\right)v-k\frac{v^2}{3}$ and it is immediately seen that if $\alpha \geq k/2$, then the assumptions of case 1) in the statement of Theorem \ref{thm-conv1} are satisfied. However, since $\lim_{s \to 0} f(s)/h'(s)=1/(2\alpha)$, the inviscid profile $\mathcal{V}_I$ exits the interval $[0, 1]$, while $\bar{v}$ is constrained between $0$ and $1$ and hence $\bar{v} \equiv 0$  on the left of the vanishing point for $\mathcal{V}_I$. In particular, $\bar{v}$ is not $C^1$ and is sharp near the value $0$, similarly to the case without convection. This outcome is reproduced in Figure \ref{fig:1}, top line. If instead $\alpha \in (-k/6, k/2)$, then $S'(1) < 0$ and case 2) of Theorem \ref{thm-conv1} occurs; moreover, it can be easily checked that $f-h'$ has a unique maximum, so that the limit profile is again sharp on the ``left'' side only (Figure \ref{fig:1}, middle line). Finally, if $\alpha \leq -k/6$, it is immediate to see that $S'(v) > 0$ for every $v \in (0, 1)$, so that the limit configuration for $v_\eps$ is fully piecewise linear, in accord with case 3) of Theorem \ref{thm-conv1} (Figure \ref{fig:1}, bottom line).  Summarizing, Theorem \ref{thm-conv1} completely characterizes $\bar{v}$ for Fisher-Burgers type equations: 
the limit profile is never regular, becoming fully sharp (near both the values $0$ and $1$) if the convection is negative and sufficiently large (in particular, \emph{concave}). Of course, for fixed $\eps > 0$ the profiles are always smooth, anyway in the pictures we can spot quite neatly the asymptotic trend stated in the theorem. 

\begin{figure}[h!]
\centering
\includegraphics[scale=0.48]{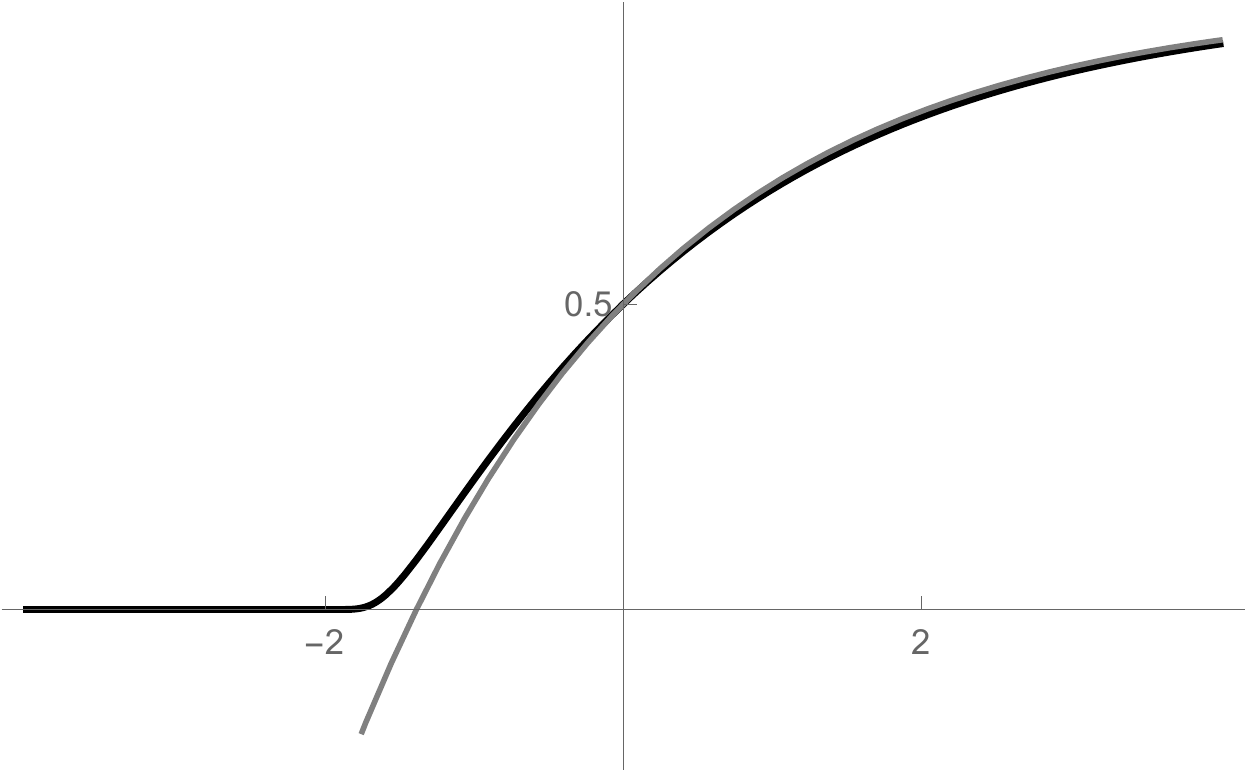}
    \hspace{1cm} \includegraphics[scale=0.48]{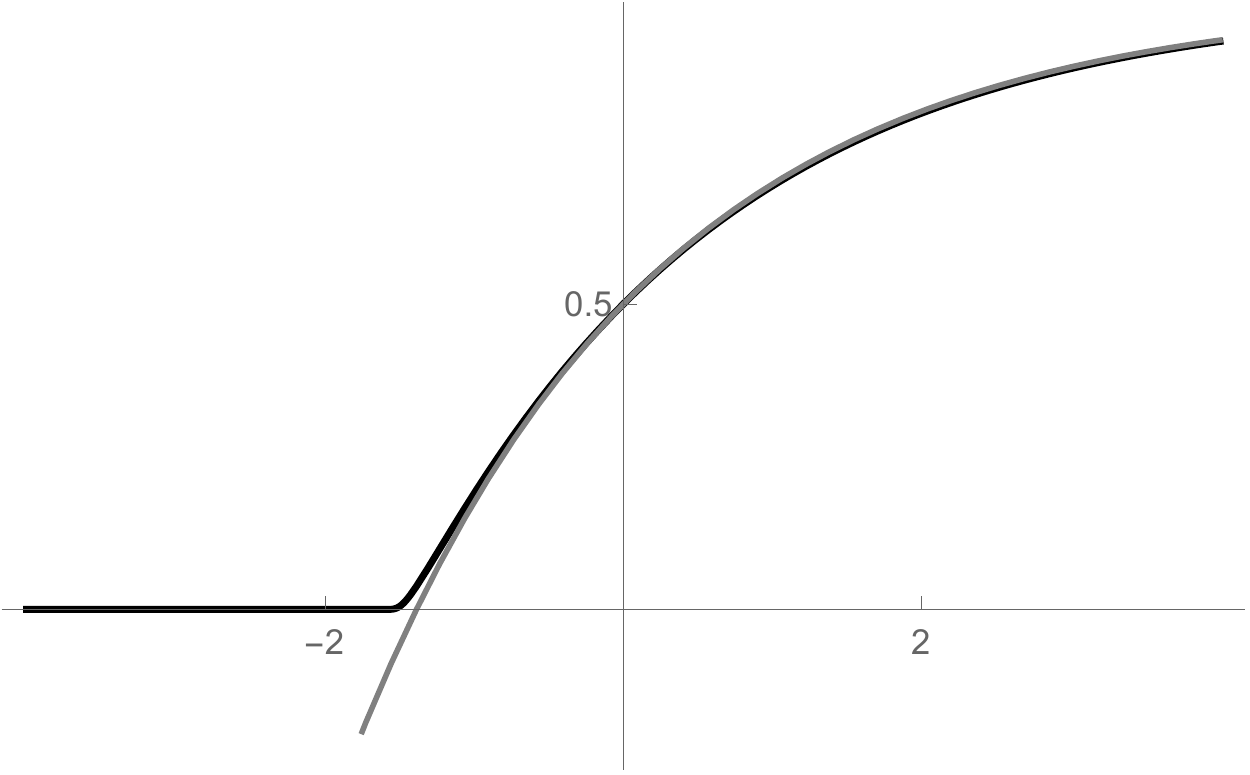}
\caption{For $k=1$ and $\alpha=1$, we depict the critical front solution $v_\eps$ of \eqref{esempio} (black) and the inviscid profile $\mathcal{V}_I$ (gray) in the cases $\eps=0.002$ (left) and $\eps=0.0002$ (right). In the former case $c_\eps^* \approx 0.07$, as in Figure \ref{fig:1}, while in the latter one we find $c_\eps^* \approx 0.024$.}\label{fig:2}
\end{figure}

In Figure \ref{fig:2}, we zoom into the above case 1), showing how the critical profile $v_\eps$ modifies its shape for smaller $\eps$; we explicitly plot the inviscid profile $\mathcal{V}_I$ (gray) to highlight how much $v_\eps$ and $\mathcal{V}_I$ become almost indistinguishable whenever strictly positive.

\begin{figure}[h!]
\centering
\includegraphics[scale=0.46]{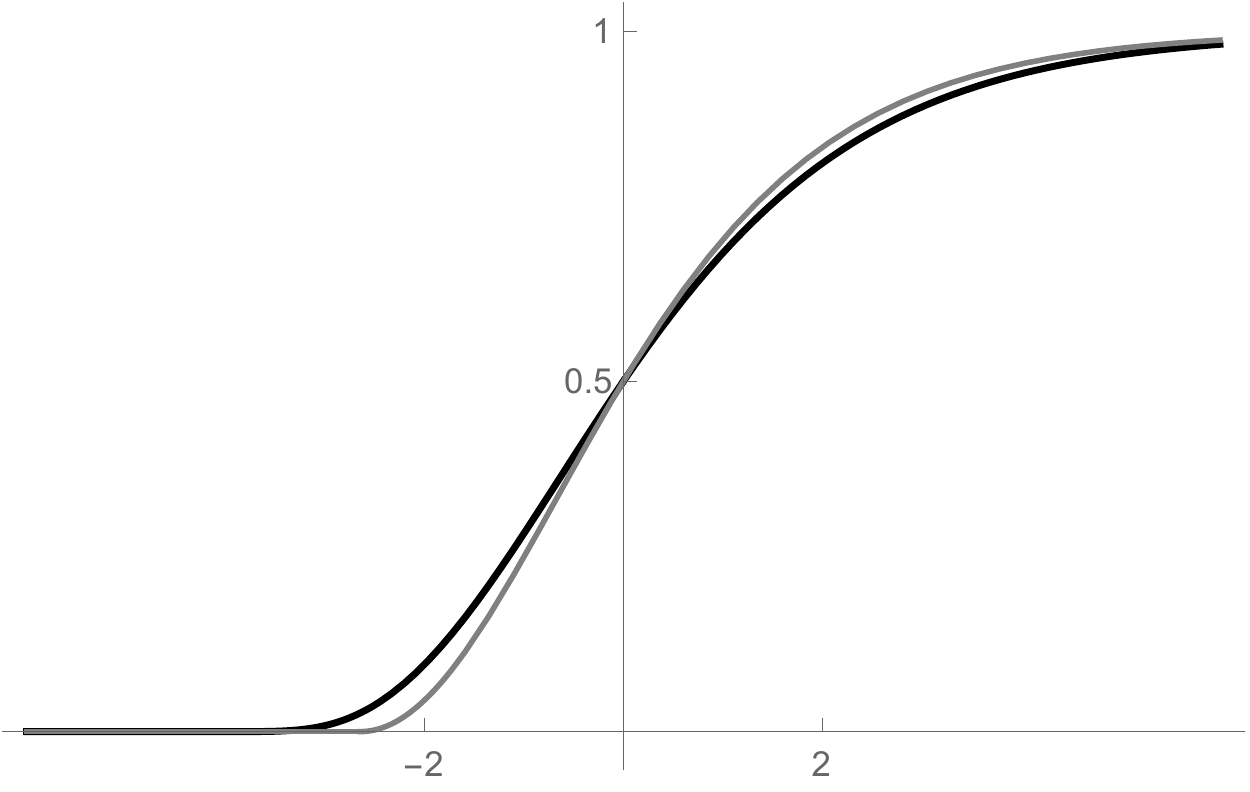}
    \hspace{1cm} \includegraphics[scale=0.46]{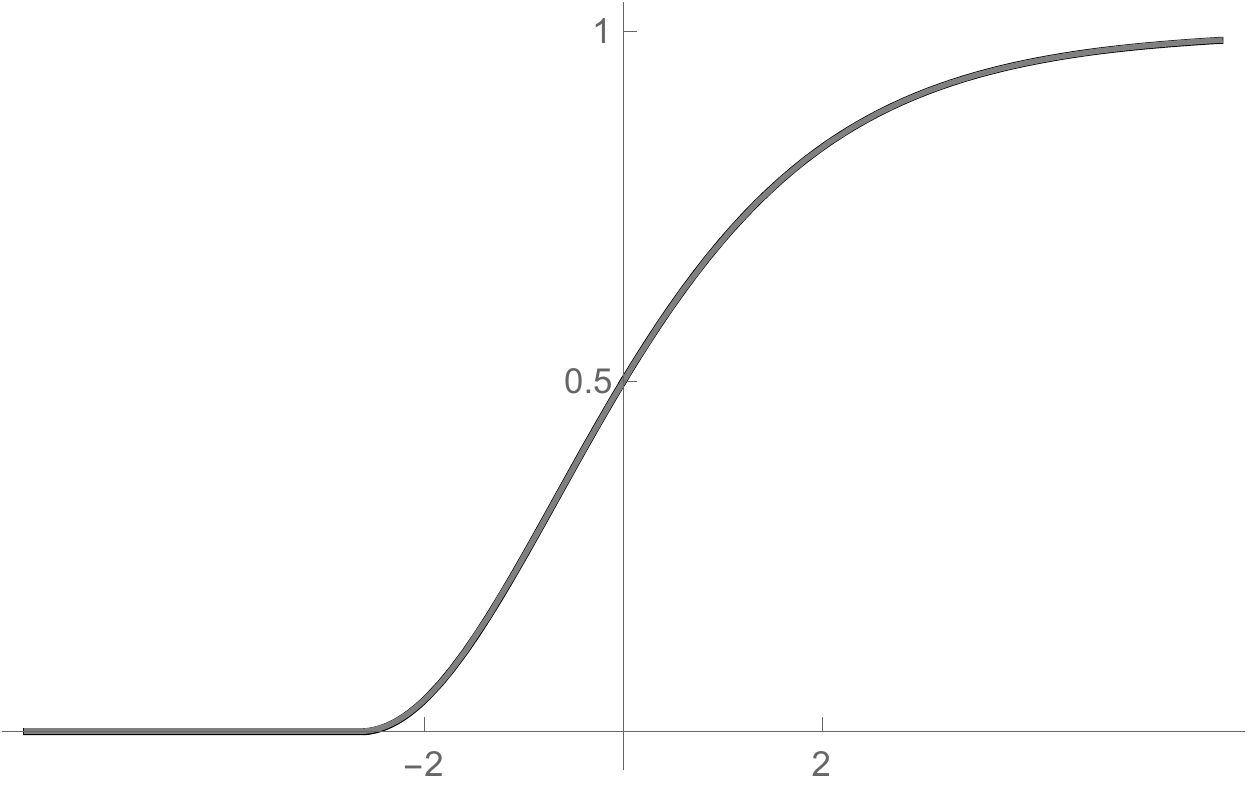}
\caption{For $h(s)=s^{3/2}$ and $f(s)=s(1-s)$, we depict the critical front profile $v_\eps$ for \eqref{MinkIntro} (black) and the inviscid profile $\mathcal{V}_I$ (gray) for $\eps=0.01$ (left) and $\eps=0.002$ (right). Here $c_\eps^* \approx 0.143$ (left), $c_\eps^* \approx 10^{-3}$ (right).}\label{fig:3}
\end{figure}

\begin{figure}[h!]
\centering
\includegraphics[scale=0.46]{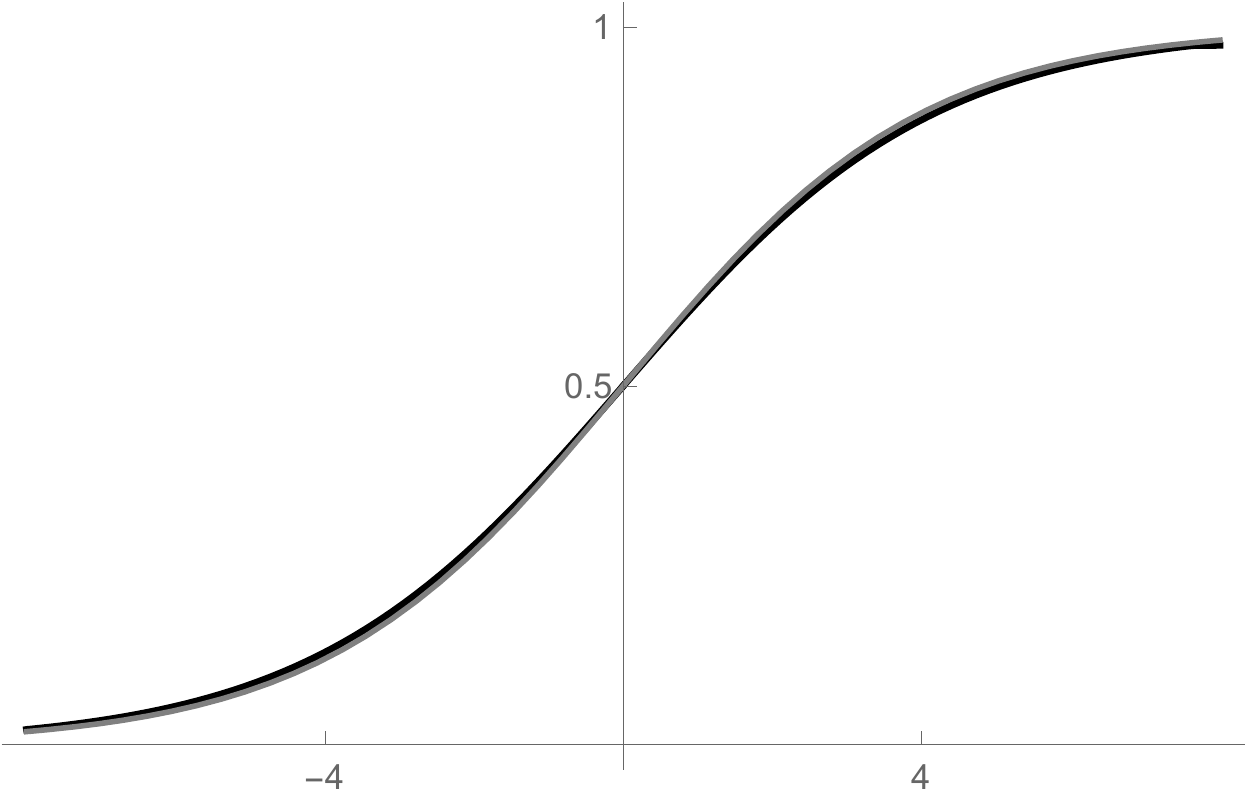}
    \hspace{1cm} \includegraphics[scale=0.46]{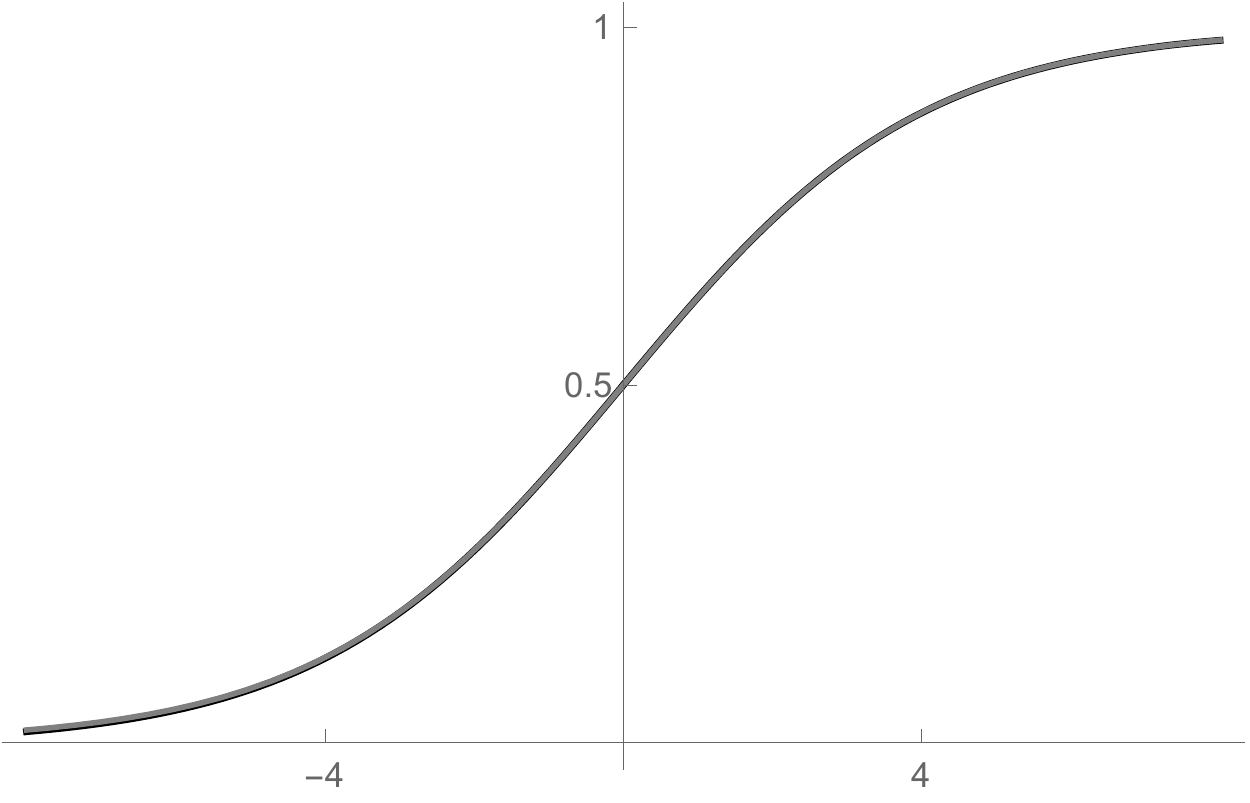}
\caption{For $h(s)=s^2$ and $f(s)=s^2(1-s)$, we depict the critical front profile $v_\eps$ for \eqref{MinkIntro} (black) and the inviscid profile $\mathcal{V}_I$ (gray) for $\eps=0.1$ (left) and $\eps=0.01$ (right). Here $c_\eps^* \approx 0.046$ (left), $c_\eps^* \approx 10^{-4}$ (right).}\label{fig:4}
\end{figure}

We finally corroborate our discussion about the importance of the value of $\lim_{s \to 0^+} f(s)/h'(s)$ in the possible regularization of the limit profile in the above case 1). In Figure \ref{fig:3}, we show the critical profile $v_\eps$ for $\eps=0.01$ (left) and $\eps=0.002$ (right) in case $f(s)=s(1-s)$ and $h(s)=s^{3/2}$; here, $\lim_{s \to 0^+} f(s)/h'(s)=0$, so that $\bar{v}$ is everywhere $C^1$ and reaches the equilibrium $0$ in finite time. On the other hand, in Figure \ref{fig:4} we show the shape of $v_\eps$ for $\eps=0.1$ (left) and $\eps=0.01$ (right) in case $f(s)=s^2(1-s)$ and $h(s)=s^2$, for which the inviscid profile $\mathcal{V}_I$ is regular and reaches the equilibrium $0$ only at $-\infty$, since the integral $\int_0^{1/2} h'(s)/f(s) \, ds$ diverges. Consequently, $\bar{v}$ is everywhere $C^1$ and $\bar{v}(z) > 0$ for every $z \in \mathbb{R}$. We notice again that $\bar{v}$ and $\mathcal{V}_I$ appear indistinguishable in the considered interval.



\bigbreak
\textbf{Acknowledgment.} This research was funded by the National Plan for Science, Technology and Innovation (MAARIFAH), King Abdulaziz City for Science and Technology, Kingdom of Saudi Arabia, Award Number (13-MAT887-02).

\end{document}